\documentclass[11pt]{article}

\usepackage{amssymb}
\usepackage{mathrsfs}

\usepackage{amsmath}

\newcommand{\funbi}[5]
{
{#1} \left\{
\begin{array}{ll}
{#2} & \mbox{{#3}} \\
{#4} & \mbox{{#5}}
\end{array}
\right.
}

\newtheorem{exam}{Example}
\newtheorem{theo}{Theorem}

\newtheorem{lemm}{Lemma}

\newtheorem{rema}{Remark}

\newcommand{\proof}{\vspace{0.00cm}{\bf Proof: }}

\title{
Inconsistency of Measure-Theoretic Probability}
\author{
\vspace{0.00cm}
Guang-Liang Li 
\footnote{Corresponding author (alternative e-mail address: guanglli@hotmail.com).
}
\hspace{2cm}
Victor O.K. Li
\\
Department of Electrical and Electronic Engineering\\
The University of Hong Kong \\
\{glli,vli\}@eee.hku.hk
}
\begin{document}
\maketitle





\begin{abstract}
We reveal a contradiction in measure-theoretic probability.
The contradiction is an ``equation'' $1/2 = 0$ with its
two sides representing probabilities. Unlike known paradoxes
in mathematics, the revealed contradiction cannot be explained away and
actually indicates that measure-theoretic probability is inconsistent.
Appearing only in the theory,
the contradiction does not exist in the physical world. So
practical applications of measure-theoretic probability will not
be affected by the inconsistency as long as ``ideal events'' in the theory
(which will never occur physically) 
are not mistaken for real events in the physical world. Nevertheless,
the inconsistency must be resolved. Constructive mathematics can avoid
such inconsistency. There is no contradiction reported in constructive
mathematics.
\end{abstract}

\section{Introduction}
\hskip\parindent
Known paradoxes in mathematics  
are either considered issues of mathematical 
logic \cite{Bena}, or treated as deducible
statements which are counterintuitive but 
true \cite{Czyz, Szek}, where ``true'' implies ``without contradiction''. 
Treating deducible statements as 
results without contradiction is not grounded
on mathematical justification. Such treatment is usually based 
on a popular belief, i.e., 
mathematics is free from contradiction as it is.
The belief is an issue of mathematical philosophy, 
which is not our concern.
In this paper, we show that a contradictory ``equation'' $1/2 = 0$ 
is deducible in measure-theoretic probability, where
the two sides of the ``equation''
represent probabilities. Unlike the paradoxes discussed in \cite{Czyz, Szek},
the contradictory ``equation'' cannot be explained away and actually indicates that
measure-theoretic probability is inconsistent.

After a brief discussion about why measure-theoretic probability is very successful and
yet inconsistent (Section \ref{sec-2}), we review some key notions relevant to revealing the inconsistency, such as
inadequacy of the $\sigma$-algebra generated by coordinate variables of stochastic
processes (Section \ref{sec-3}), weak convergence and tightness of probability measures on the $\sigma$-algebra
${\mathcal B}({\mathbb R})$ of subsets of
the real line ${\mathbb R}$ (Section \ref{sec-4}), weak convergence of probability measures on the $\sigma$-algebra 
${\mathcal B}(\overline{\mathbb R})$ of subsets of
the extended real line $\overline{\mathbb R}$ (Section \ref{sec-5}), and the difference between probabilities and
numbers in the closed unit interval $[0, 1]$ (Section \ref{sec-6}), in particular, 
the difference between limiting probabilities and limits of sequences 
in $[0, 1]$ (Section \ref{sec-7}).
Based on the above preparation, we demonstrate
the deducibility of the contradictory ``equation'' $1/2 = 0$ (Section \ref{sec-8}). 
We conclude in Section \ref{sec-9}.

\section{Measure-Theoretic Probability: Success and Inconsistency}
\label{sec-2}
\hskip\parindent
For measure-theoretic probability, 
the ultimate criterion of success 
is whether the probability of a real event (i.e., an event in the physical world)
predicted by the theory
can be verified empirically. By this criterion, measure-theoretic
probability is very successful. And yet,
measure-theoretic probability is 
inconsistent, as we shall demonstrate. 
How can measure-theoretic probability 
be both successful and inconsistent? A brief answer is as follows.
As a consequence of ``ideal events'',
the inconsistency of measure-theoretic probability 
appears only in the theory, but 
not in the physical world. 
Expressed as an
infinite union or intersection of different
events of positive probability, an 
``ideal event'' is a result of infinitely repeated or
infinitely fine operations.

\begin{exam}
\label{ex-ideal1}
{\em
An absolutely accurate measurement of a physical
quantity with its relative magnitude represented by an irrational number
is an ideal event. Corresponding to picking up a point from a set with the
cardinality of the
continuum, such measurement requires
infinitely fine operations,
and hence is not physically realizable. 
$\Box$
}
\end{exam}

The criterion of success of measure-theoretic probability 
implies a constraint: Since ``ideal events'' do not exist physically, 
or since it is impossible to verify the
existence of ``ideal events'' in the physical world,
physical existence of theoretical properties derived by
manipulating
probabilities of ``ideal events''
should not be inferred from experiments.
This constraint serves as
a ``firewall'' of measure-theoretic probability. When the constraint is satisfied,
it is safe to apply measure-theoretic probability to solve problems in the real world.
In practice, however,
violating the constraint can cause
confusion of real and ideal events.
A typical violation is to use 
a finite set of experimental data to ``prove'' or
``verify''
physical existence of theoretical
properties obtained by calculating probabilities of ``ideal events''.

\section{Inadequacy of $\sigma$-Algebra Generated by Coordinate Variables}
\label{sec-3}
\hskip\parindent
A real-valued random variable is usually described by the 
probability measure (i.e., the distribution) induced by the random variable on the real line ${\mathbb R}$.
In general, a stochastic process is usually described by the probability measures 
(i.e., the finite-dimensional distributions) induced by the process in $k$-dimensional
Euclidean spaces ${\mathbb R}^k, k = 1, 2, \cdots$. 
The construction of a process in measure-theoretic probability begins with a given
system of finite-dimensional distributions. However, the system of
finite-dimensional distributions may not completely capture the
properties of a process \cite{Bill}.

A stochastic process is a collection of real-valued random variables labeled by an
index variable taking values from a set $T$. The points in the set $T$ usually represent time.
We consider only $T=\{1, 2, \cdots\}$. Denote by
${\mathbb R}^T$ the set of all real sequences.
For each $n\in T$, define a mapping $W_n: {\mathbb R}^T\to {\mathbb R}$ by
\[
W_n(\omega) = \omega(n) = \omega_n.
\]
For a given $n\in T$, the mapping $W_n$ (also referred to as a coordinate function)
maps a real sequence $\omega\in {\mathbb R}^T$ to the $n$th coordinate $\omega(n) = \omega_n$ of $\omega$. 
We may also write $W_n(\omega)$ as $W(n, \omega)$. For a fixed $\omega$, 
$W(\cdot, \omega) = \omega(\cdot)$ is just $\omega$ itself.
Denote by ${\mathcal B}({\mathbb R}^T)$ 
the $\sigma$-algebra (also referred to as $\sigma$-field)
generated by all the coordinate functions $W_n, n\in T$:
${\mathcal B}({\mathbb R}^T) = \sigma\{W_n: n\in T\}$, which is the intersection of all the $\sigma$-algebras
containing the sets of the form 
\[
\{\omega\in {\mathbb R}^T: W_n(\omega)\in H\} = \{\omega\in {\mathbb R}^T: \omega_n\in H\}
\]
for $n\in T$ and $H\in {\mathcal B}({\mathbb R})$. In other words,
${\mathcal B}({\mathbb R}^T)$ is the ``smallest'' class of the sets which is closed under the finite and countable
set-theoretic operations.

Let ${\mathbb R}^T$ be a sample space represented by $\Omega$, and
let a probability measure be
defined on ${\mathcal B}({\mathbb R}^T)$. 
Now $\omega\in \Omega = {\mathbb R}^T$ represents a sample point, and $W_n$ can be interpreted as
the $n$th coordinate variable $U_n$ of a random sequence
$(U_n)_{n\geq 1}$. 
A realization of $(U_n)_{n\geq 1}$ at a given sample point $\omega$ is
a real sequence $(U_n(\omega))_{n\geq 1}$, which can also be expressed as
$U(\cdot, \omega)$ and is just $\omega$ itself. 

If $(U_n)_{n\geq 1}$ is a non-decreasing random sequence, 
$U_n(\omega)$ may tend to infinity as $n\to\infty$ at
some $\omega$. Since
$\{\infty\}\not\in {\mathcal B}({\mathbb R}^T)$, 
a probability measure
on the measurable space 
$({\mathbb R}^T, {\mathcal B}({\mathbb R}^T))$ 
is not useful
to describe the limiting behavior of $(U_n)_{n\geq 1}$.
Clearly, to define the limiting random variable of 
$(U_n)_{n\geq 1}$ on a probability space, the
limiting random variable should be measurable with respect
to the corresponding probability measure,
which requires the values of the limiting random variable to
form a measurable set on the probability space.
Fortunately, 
the probability measures 
$\rho_n$ induced by $U_n$ constitute a sequence $(\rho_n)_{n\geq 1}$
such that,
corresponding to $(U_n)_{n\geq 1}$ and $(\rho_n)_{n\geq 1}$, 
a limiting random variable $U$ and a
limiting probability measure $\rho$ induced by $U$ exist on the
extended real line $\overline{\mathbb R}$ and the corresponding $\sigma$-algebra
${\mathcal B}(\overline{\mathbb R})$, respectively.
For a 
convergent sequence of probability
measures, the corresponding  
random variables can be on different probability spaces \cite{Bill2}. Consequently,
convergence of probability measures does not require
the involved random variables to form a convergent
sequence on the same probability space. Convergence of probability measures on ${\mathcal B}(\overline{\mathbb R})$ is a generalization
of convergence of probability measures on 
${\mathcal B}({\mathbb R})$. If a sequence of probability measures converges on ${\mathcal B}({\mathbb R})$, the sequence is necessarily
tight on ${\mathcal B}({\mathbb R})$. We shall review these notions
in the next two sections.

\section{Weak Convergence and Tightness of Probability Measures on 
${\mathcal B}({\mathbb R})$}
\label{sec-4}
\hskip\parindent
Let $(\rho_n)_{n\geq 1}$ be a sequence of probability measures on ${\mathcal B}({\mathbb R})$. The probability measures $\rho_n$ are distributions induced by 
real-valued random variables $U_n$, such that
$U_n$ converge in distribution to a limiting random variable $U$
as $n\to\infty$ (notation $U_n\Rightarrow U$).
Equivalently,
the probability measures $\rho_n$
converge weakly to
a limiting probability measure $\rho$
induced by $U$ (notation $\rho_n\Rightarrow\rho$). 
If $(U_n)_{n\geq 1}$ is non-decreasing,
$U$ may take $\infty$
as its value with positive probability, i.e., $\rho$
may assign positive probability to
$\{\infty\}$. 
One way to define convergence in distribution of
a sequence of real random variables is by convergence of the corresponding
distribution functions at each continuous point of the limiting distribution function,
which may not necessarily be a distribution function of a real random variable \cite{Loev}.
This definition has some disadvantages, for example, see \cite{Poll}.
 
An equivalent definition of weak convergence of 
probability measures (i.e., convergence in distribution of 
the corresponding random variables)
is as follows \cite{Bill2,Poll}.
Suppose ${\mathbb R}$ is the common range for all the random variables involved.
The sequence $(\rho_n)_{n\geq 1}$ of probability measures induced by the
random variables converges weakly to a
probability measure $\rho$, if and only if for each $f\in {\mathbb C}_b({\mathbb R})$,
\[
\int_{\mathbb R}fd\rho_n\to \int_{\mathbb R}fd\rho,\; n\to\infty
\]
where ${\mathbb C}_b({\mathbb R})$ is the class of bounded, continuous functions on ${\mathbb R}$.
Clearly, such convergence is not only a property of
the probability measures $\rho_n$ or corresponding random variables $U_n$, but also relevant to the common range space
of the involved
random variables on which weak convergence is defined. 
If the common range space is defined to be ${\mathbb R}$, and if the limiting random variable is not
real-valued, then weak convergence or convergence in distribution as defined above does not hold, even if
a limiting random variable and the corresponding probability measure exist on $\overline{\mathbb R}$ and 
${\mathcal B}(\overline{\mathbb R})$, respectively.

\begin{exam}
\label{ex-Dirac}
{\em
Let $(U_n)_{n\geq 1}$ be a sequence of degenerate random variables such that
$U_n = n$. The distribution, i.e., the probability measure $\rho_n$ induced by $U_n$ is
the Dirac measure concentrated at $\{n\}$. As $n\to\infty$, 
$U_n\Rightarrow U =\infty\in\overline{\mathbb R}$ and
$\rho_n\Rightarrow\rho$ where $\rho$ is the Dirac measure concentrated on 
$\{\infty\}\in{\mathcal B}(\overline{\mathbb R})$. But $(U_n)_{n\geq 1}$ does not converge on ${\mathbb R}$, and
$(\rho_n)_{n\geq 1}$ does not converge weakly on ${\mathcal B}({\mathbb R})$. 
$\Box$
}
\end{exam}

On the one hand,
as illustrated by the above example,
for a sequence of probability measures, mass may ``escape to infinity'' as $n\to\infty$,
and meanwhile, in an extended sense,
the notion of
limiting probability measure remains meaningful. 
In particular, for a non-decreasing random sequence such as $(U_n)_{n\geq 1}$ in 
Example \ref{ex-Dirac}, a unique
limiting random variable exists, which may not necessarily be
real-valued. By Helly's selection theorem (e.g., see \cite{Bill}), the limiting random variable has 
unique distribution function and probability measure, such that the limiting distribution function may
have a lower bound larger than 0 or an upper bound less than 1, and the limiting probability measure $\rho$ is on
${\mathcal B}(\overline{\mathbb R})$ and may be a subprobability measure when restricted to ${\mathcal B}({\mathbb R})$,
i.e., $\rho({\mathbb R}) < 1$. The cause of $\rho({\mathbb R}) < 1$ is mass ``escaping to infinity''.

On the other hand, probability measures defined on ${\mathcal B}({\mathbb R})$ are particularly
important for measure-theoretic probability and its applications.
To prevent mass from ``escaping to infinity'',
a condition must be imposed on a sequence of probability measures before we consider
whether the sequence converges to a limiting probability measure on ${\mathcal B}({\mathbb R})$.

Tightness of probability measures on ${\mathcal B}({\mathbb R})$ is the condition to prevent
mass from ``escaping to infinity''.
If a sequence of probability measures
$(\rho_n)_{n\geq 1}$ is tight on
$({\mathbb R}, {\mathcal B}({\mathbb R}))$, then
for each $\epsilon > 0$, there is an interval $I_{\epsilon}\subset {\mathbb R}$,
such that $\sup_n\rho_n(I_{\epsilon}^c) < \epsilon$ where 
$I_{\epsilon}^c = {\mathbb R}\setminus I_{\epsilon}$. If such a tight sequence
of probability measures $(\rho_n)_{n\geq 1}$ converges weakly to $\rho$,
then the weak limit
$\rho$ is also a probability measure on $({\mathbb R}, {\mathcal B}({\mathbb R}))$.
Given $f\in {\mathbb C}_b({\mathbb R})$, even if $f(x)$ does not converge as $x\to \pm\infty$,
so long as $(\rho_n)_{n\geq 1}$ is tight on $({\mathbb R}, {\mathcal B}({\mathbb R}))$
and $\rho_n\Rightarrow\rho$, we have
$\rho(\{\pm\infty\}) = 0$ and 
$\int_{{\mathbb R}} fd\rho$ exists. 
Suppose $\rho_n\Rightarrow\rho$. If not only all $\rho_n$ but also $\rho$ are
required to be probability measures on $({\mathbb R}, {\mathcal B}({\mathbb R}))$,
then $(\rho_n)_{n\geq 1}$ is necessarily tight on $({\mathbb R}, {\mathcal B}({\mathbb R}))$.

\begin{exam}
\label{ex-tight}
{\em
Let $\lambda, \lambda_n, n = 1, 2, \cdots$ 
be probability measures on 
$({\mathbb R}, {\mathcal B}({\mathbb R}))$, 
such that $\lambda, \lambda_1, \lambda_2, \cdots$ are
distributions of i.i.d. random variables $X, X_1, X_2, \cdots$, 
respectively. The probability mass functions of $X$ and $X_n$ are the same:

\[
\lambda(\{j\}) =\lambda_n(\{j\}) =  \left\{
\begin{array}{llr}
1/2, & j = 0\\
1/2, & j = 1\\
0,   & \mbox{otherwise.} 
\end{array}\right.
\]
If we define $\lambda(\{\pm\infty\})=\lambda_n(\{\pm\infty\})=0$, then $\lambda, \lambda_n, n = 1, 2, \cdots$ 
are also probability measures on 
$(\overline{{\mathbb R}}, {\mathcal B}(\overline{\mathbb R}))$ trivially. 
The sequence $(\lambda_n)_{n\geq 1}$ is tight on $({\mathbb R}, {\mathcal B}({\mathbb R}))$, and
converges weakly to $\lambda$. Equivalently, $X_n$ converge in distribution to $X$.
$\Box$
}
\end{exam}

\begin{exam}
\label{ex-notight}
{\em
Consider a random sequence
$(Z_n)_{n\geq 1}$. 
Denote by $\mu_n$ the probability measure (i.e., the distribution) induced by $Z_n$, and let 
the probability mass function of $Z_n$ be
\begin{equation}
\label{eq-1}
\mu_n(\{j\}) = \left\{
\begin{array}{llr}
2^{-(n - j + 1)}, & 1\leq j < n\\
2^{-1} + 2^{-n},  & j = n\\
0,                & \mbox{otherwise.} 
\end{array}\right.
\end{equation}
For each $n$, $Z_{n+1}$ takes larger values with larger probabilities.
The values of $Z_{n+1}$ include $n+1$ and all the values of $Z_n$.

The sequence $(\mu_n)_{n\geq 1}$ is not tight on $({\mathbb R}, {\mathcal B}({\mathbb R}))$.
To see this,
consider $\epsilon = 2^{-1}$. By (\ref{eq-1}), 
for all $n$, $\mu_n(\{n\}) = 2^{-1} + 2^{-n} > 2^{-1}$.
In particular, for all real number $b > 0$ and all $n > N_b$, 
$\mu_n(\{n\}) > 2^{-1}$, where $N_b= \lceil b\rceil$, and
$\lceil b\rceil$ is the smallest integer larger than or equal to $b$.
Consequently, for each interval $I_b=(-\infty, b)$, 
$\sup_n\mu_n(I_b^c) > 2^{-1}$.
The probability mass function (\ref{eq-1}) 
explains not only why $(\mu_n)_{n\geq 1}$ fails to be tight on $({\mathbb R}, {\mathcal B}({\mathbb R}))$
but also where all the mass goes as $n\to\infty$.
For a fixed $j$, $\lim_n\mu_n(\{j\})=0$. 
To see what happens when $j$ increases with $n$, 
write $j = n - k$, where $0\leq k < n$ is fixed.
We then have
$\mu_n(\{n-k\})=1/2^{k+1}$ for $k > 0$ and 
$\mu_n(\{n\})=2^{-1} + 2^{-n}$.
Since $(n-k)\to\infty$ as $n\to\infty$,  
$1/2^{k+1}, k = 1, 2, \cdots$ constitute one half
of the mass that
``escapes to infinity''. 
The other half of the mass at $\{\infty\}$ in the
limit as $n\to\infty$
comes from $\mu_n(\{n\})$.
Let a vertical line segment of
length $\mu_n(\{j\})$ at $j$ represent
the mass assigned to $\{j\}$ by $\mu_n$.
As $n\to\infty$, 
the length of the segment at a fixed $j$ diminishes towards 0, the segment at $j = n - k$ with a fixed $0< k < n$
maintains an unchanged length $2^{-(k + 1)}$, and the length of the segment at $n$ is $2^{-1} + 2^{-n}$,
which decreases gradually to 1/2. 
Figure 1 illustrates the case for $k > 0$ with $n = 2$ and 3. 
The total ``escaped'' mass then is
$\sum_{k\geq 0}1/2^{k+1} = 1$. 
Actually, by Helly's theorem, each subsequence of $(\mu_n)_{n\geq 1}$ converges
to a subprobability measure $\mu$ on $({\mathbb R}, {\mathcal B}({\mathbb R}))$ such that $\mu({\mathbb R}) = 0$,
for $\mu({\mathbb R}) > 0$ contradicts (\ref{eq-1}). So $\mu_n\Rightarrow\mu$ where 
$\mu$ is the Dirac measure on 
$(\overline{\mathbb R}, {\mathcal B}(\overline{\mathbb R}))$
with $\mu(\{\infty\}) = 1$. 
Equivalently, $Z_n\Rightarrow Z$, where $Z$ is the
extended random variable corresponding to $\mu$.
When restricted to $({\mathbb R}, {\mathcal B}({\mathbb R}))$, $\mu$ is
a zero measure. 
$\Box$
}
\end{exam}

As we shall see later, 
$(\lambda_n)_{n\geq 1}$ in Example \ref{ex-tight} has an intrinsic connection to 
$(\mu_n)_{n\geq 1}$ in Example \ref{ex-notight}, though the former is tight 
on $({\mathbb R}, {\mathcal B}({\mathbb R}))$ and the
latter is not.

\begin{figure*}
\setlength{\unitlength}{0.4mm}
\begin{picture}(100,110)

\put(20,60){\vector(1,0){100}}
\put(62,60){\line(0,1){26}}
\put(62,93){\makebox(0,0){$1/4$}}
\put(50,52){\makebox(0,0){$Z_2 = 1$}}
\put(128,60){\makebox(0,0){$x$}}

\put(160,60){\vector(1,0){100}}
\put(189,60){\line(0,1){13}}
\put(230,60){\line(0,1){26}}
\put(189,80){\makebox(0,0){$1/8$}}
\put(230,95){\makebox(0,0){$1/4$}}
\put(177,52){\makebox(0,0){$Z_3=1$}}
\put(219,52){\makebox(0,0){$Z_3=2$}}
\put(268,60){\makebox(0,0){$x$}}

\put(295,60){\vector(1,0){100}}
\put(350,52){\makebox(0,0){$Z<\infty$}}
\put(403,60){\makebox(0,0){$x$}}

\put(210,25){\makebox(0,0){Figure 1. 
For any fixed $k\geq 0$, mass at $\{n-k\}$ ``escapes
to infinity'' as $n\to\infty$, so 
$\mu({\mathbb R})= 0$.}}
\end{picture}
\end{figure*}

\section{Weak Convergence of Probability Measures on 
${\mathcal B}(\overline{\mathbb R})$}
\label{sec-5}
\hskip\parindent
The definition of $\rho_n\Rightarrow\rho$ can be generalized 
such that, for each $n$, $\rho_n$ is still a probability measure on 
$({\mathbb R}, {\mathcal B}({\mathbb R}))$, but $(\rho_n)_{n\geq 1}$ may not
be tight on $({\mathbb R}, {\mathcal B}({\mathbb R}))$, and $\rho$ is a probability measure
on $(\overline{\mathbb R}, {\mathcal B}(\overline{\mathbb R}))$ with $\rho(\{\pm\infty\}) > 0$ (see Example \ref{ex-Dirac}).
In the above generalized definition of $\rho_n\Rightarrow\rho$, ${\mathbb C}_b({\mathbb R})$ in the
original definition is replaced by ${\mathbb C}(\overline{\mathbb R})$, the class of continuous functions
on $\overline{\mathbb R}$. The reason for such replacement is shown
by the example below. Let ${\bf 1}_A$ represent the indicator function, i.e.,
\[
\funbi{{\bf 1}_A(x) =}{1,}{$x\in A$}
{0,}{otherwise.}
\]
\begin{exam}
\label{ex-Dirac2}
{\em
Consider $(\rho_n)_{n\geq 1}$ such that $\rho_n$ is the Dirac measure
concentrated at $\{n\}$ as in Example \ref{ex-Dirac}. The distribution function corresponding to $\rho_n$ is
${\bf 1}_{[n, \infty)}(\cdot)$ with ${\bf 1}_{[n, \infty)}(x)\to 0$ as $n\to\infty$
for $\lim_n[n, \infty)=\cap_n[n, \infty)=\emptyset$. By Helly's selection theorem and
the generalized definition of $\rho_n\Rightarrow\rho$, $\rho$ is the Dirac
measure concentrated at $\{\infty\}$. However, for some $f\in {\mathbb C}_b({\mathbb R})$ such as $f =\sin$,
$\int_{{\mathbb R}} fd\rho_n = f(n) = \sin(n)$, which does not converge 
as $n\to\infty$. 
$\Box$ 
}
\end{exam}

Let $\overline{\mathbb R}$ be equipped with a topology induced by a fixed metric.
A function from one metric space to another is called a homeomorphism if the function
is continuous, one-to-one, onto, and its inverse function is also continuous.
If such a function exists, then the two metric spaces are said to be homeomorphic.
With a homeomorphism from $\overline{\mathbb R}$ onto
$[0, 1]$, ${\mathbb C}(\overline{\mathbb R})$ can be identified with ${\mathbb C}([0, 1])$,
the space of continuous functions on $[0, 1]$. 
\begin{exam}
\label{ex-homeo}
{\em
A homeomorphism $h: \overline{\mathbb R}\to[0, 1]$ from $\overline{\mathbb R}$ onto
$[0, 1]$ is
\[
h(x) =  \left\{
\begin{array}{llr}
0, & x = -\infty\\
1-\pi^{-1}\cot^{-1}(x), & x\in{\mathbb R}\\
1,   & x=\infty.
\end{array}\right.
\]
The inverse of $h$ is $h^{-1}: [0, 1]\to\overline{\mathbb R}$ given by
\[
h^{-1}(y) =  \left\{
\begin{array}{llr}
-\infty, & y = 0\\
\cot(\pi(1-y)), & y\in(0, 1)\\
\infty,   & y=1.
\end{array}\right.
\]
A usual metric on $\overline{\mathbb R}$ is $d(x, y) = |\tan^{-1}(x) - \tan^{-1}(y)|$.
$\Box$
}
\end{exam}

Each $f\in {\mathbb C}([0, 1])$ has limits
at the boundary points.
Weak convergence of probability measures
on ${\mathcal B}(\overline{\mathbb R})$ can then be identified with
weak convergence of probability measures on ${\mathcal B}([0, 1])$
(the $\sigma$-algebra of subsets of $[0, 1]$), and hence
a limiting probability measure 
on ${\mathcal B}([0, 1])$ which assigns mass to a boundary point of $[0, 1]$
is identified with the corresponding limiting probability measure
on ${\mathcal B}(\overline{\mathbb R})$ which assigns the same amount of mass to
$\{\infty\}$. Since the probability measures on ${\mathcal B}([0, 1])$
may not necessarily be tight on ${\mathcal B}((0, 1))$ (the $\sigma$-algebra of subsets of
$(0, 1)$),
weak convergence of such probability measures on ${\mathcal B}([0, 1])$
does not necessarily imply tightness of the corresponding
probability measures on ${\mathcal B}({\mathbb R})$.
\begin{exam}
\label{ex-homeo2}
{\em
Consider $(Z_n)_{n\geq 1}$ in Example \ref{ex-notight} again. Write
$Z'_n = h(Z_n)$ where $h$ is given in Example \ref{ex-homeo}. 
Denote by $\mu'_n$ the probability measure induced by $Z'_n$.
As $n\to\infty$,
$Z_n\Rightarrow Z$. Accordingly $Z'_n\Rightarrow Z'$. Let $\mu'$ be the probability measure
induced by $Z'$. Now $\mu'_n\Rightarrow\mu'$ on
$([0, 1], {\mathcal B}([0, 1]))$ corresponds to 
$\mu_n\Rightarrow\mu$ on
$(\overline{\mathbb R}, {\mathcal B}(\overline{\mathbb R}))$.
Consequently,
$\mu'(\{1\}) = 1$ and
$\mu'((0, 1)) = 0$ correspond to 
$\mu(\{\infty\}) = 1$ and $\mu({\mathbb R}) = 0$, respectively.
Although $\overline{\mathbb R}$ and $[0, 1]$ are homeomorphic,
$(\mu_n)_{n\geq 1}$ is still not tight on $({\mathbb R}, {\mathcal B}({\mathbb R}))$. Clearly, the correspondence between 
$Z'=1$ and $Z=\infty$ does not imply that $Z'$
and $Z$ are equivalent in magnitude. The trivial tightness of $(\mu'_n)_{n\geq 1}$
on $([0, 1], {\mathcal B}([0, 1])$ will not change the fact
that $(\mu_n)_{n\geq 1}$ is not tight on $({\mathbb R}, {\mathcal B}({\mathbb R}))$. 
$\Box$
}
\end{exam}

\section{Probabilities and Real Numbers in $[0, 1]$}
\label{sec-6}
\hskip\parindent
The notion of real numbers and the notion of probability can 
both be formulated axiomatically. However, the axioms
for real numbers cannot formulate the
notion of probability. As a set function in measure-theoretic probability
formulated by the axioms of probability,
a probability measure assigns a numerical probability 
(also referred to as probability mass) represented by a number in the closed unit interval $[0, 1]$ 
to an event characterized by a measurable set on a probability space.
The numerical probability 
quantifies the possibility of the event.  In other words,
associated with a measurable set on a
specified probability space,
a (numerical) probability describes quantitatively how likely
an event will occur. 
Real numbers without the constraint imposed by the axioms of probability
do not have such unique property of probabilities.
Therefore, a number in $[0, 1]$ is not necessarily
a probability of
an event on a given
probability space.
Probabilities must comply with the constraint imposed by
the axioms of probability.
Numbers in $[0, 1]$ are not necessarily subject to
such constraint. Consequently,
although the amount or relative magnitude of a 
probability is represented by a real
number in $[0, 1]$,
probabilities cannot be treated as such numbers.

Real numbers are an ordered field. In such a field,
comparison of real numbers is defined by the axioms of order.
``Possibility'' is a quantifiable attribute of an event.
Probabilities are compared as quantified possibilities.
Since the amount of a probability is represented by a real number
in $[0, 1]$, comparing probabilities is equivalent to
comparing their amounts. 
In contrast, if a real number in $[0, 1]$ fails
to comply with the axioms of probability, then the number
is not a quantified possibility of an event on the
underlying probability space, and hence is not a probability.
Comparing such a number with a probability is meaningless. 
Such comparison is not defined
by the axioms for either probability or real numbers.  
Actually, relations on ${\mathbb R}$ 
such as $<, >, =, \leq$ and $\geq$
are all defined for and only for real numbers, where real numbers
are merely ``points''
in the ordered field. No other
interpretation is attached to such abstract points.
\begin{exam}
\label{ex-unitpro}
{\em
Unit probability is a unique quantity for
indicating certainty or almost sure occurrence of
an event.
In contrast, the real number ``1''
is merely a point 
in an absolute space ${\mathbb R}$ formulated as a set based on the axioms of real numbers without any other interpretation. 
``An event occurs with probability 1'' is certainly
a meaningful statement. However, if ``with probability 1'' is replaced by
``with the real number 1'', the statement
``an event occurs with the real number 1'' is meaningless.
Unit probability and the real number ``1'' are of the same amount, but
not of the same kind, and hence
not comparable. $\Box$
}
\end{exam}

\begin{exam}
\label{ex-info}
{\em
From information theory we can also see
the difference between probabilities and
numbers in $[0, 1]$.
As shown in information
theory, the notion of probability
is closely related to the notion of ``information''. 
Without the constraint imposed by the axioms of 
probability,
numbers in $[0, 1]$ have no such 
connection to ``information''.
$\Box$
}

\end{exam}
 
Nevertheless, since the amount of a probability is
represented by a real number in $[0, 1]$,
according to the axioms of real numbers,
we can still compare the amount of a probability with a real number.
However, when making such comparison,
we must keep in mind that a probability and its amount
(i.e., its relative magnitude represented by a real number)
are different.

An operation defined on a set is said to be
``closed'' with respect to
a given property of some or all elements in the
set, if every result of the operation is an
element of the set and possesses the property.
The values of probabilities are represented by numbers in $[0, 1]$. However,
although operations defined for numbers in $[0, 1]$ may be
closed with respect to numbers, i.e., results of
such operations are also numbers, the operations may not
necessarily be closed with respect to probabilities.
Any closed operation defined for probabilities
must obey the axioms of probability.
This requirement will not change the closeness of an operation
defined for numbers. When a closed operation defined for
numbers violates the axioms of probability,
the result of the operation is not a probability, but is
still a number.
\begin{exam}
\label{ex-4}
{\em
A probability space has four elementary events
represented by $\omega_i$ with probabilities
\[
P(\{\omega_i\}) = 1/8, i = 1, 2, 3,\; \mbox{and}\; P(\{\omega_4\}) = 5/8.
\]
For $A = \{\omega_1, \omega_2\}$ and $B = \{\omega_2, \omega_3\}$,
the condition for $P(A)$ and $P(B)$ to be additive is
violated by $A\cap B\not=\emptyset$. Consequently
$P(A) + P(B) = 1/2$ is not a probability of any event on the
probability space. Nevertheless, as the sum of two numbers, $P(A) + P(B)$
is still a number.   
$\Box$ 
}
\end{exam}

Let $\rho$ and
$\rho_n, n = 1, 2, \cdots$ be probability measures on a common measurable space
$(\overline{\mathbb R}, {\mathcal B}(\overline{\mathbb R}))$, such that
$\rho_n\Rightarrow\rho$, i.e., $\rho_n$ converges weakly to $\rho$ as $n\to\infty$,
where $\rho(\{-\infty\}) = 0$ and 
$\rho_n({\mathbb R}) = 1$ for each $n$.
Let $\alpha$ and $\alpha_n, n = 1, 2, \cdots$ be measurable sets in
${\mathcal B}(\overline{\mathbb R})$, such that $\alpha = \lim_n\alpha_n, n\to\infty$, 
$\alpha_n \in {\mathcal B}({\mathbb R})
\subset {\mathcal B}(\overline{\mathbb R})$, and
$\rho_n(\alpha_n) < 1$ for all $n$. 
Suppose $\lim_n\rho_n(\alpha_n)$ exists.
If $\alpha = {\mathbb R}$, and if $\lim_n\rho_n(\alpha_n)$ is treated
as the limit of a real sequence, i.e., if $\lim_n\rho_n(\alpha_n)$ is
interpreted as $\lim_n[\rho_n(\alpha_n)]$, 
then the result of such treatment may not necessarily
be a probability. 
\begin{exam}
\label{ex-limp}
{\em
Let $\rho_n$ be $\mu_n$ as given in Example \ref{ex-notight}.
Suppose $\alpha_n = I(x_n) = (-\infty, x_n]$. Consider $x_n = n-1$.
\[
\rho_n(\alpha_n) = \mu_n(I(n-1))
=\sum_{j=1}^{n-1}\mu_n(\{j\}) = \sum_{j=1}^{n-1}\left(\frac{1}{2^{n-j+1}}\right)
= \frac{1}{2} - \frac{1}{2^n}.
\]
Suppose $\lim_n\mu_n(I(n-1))$ is treated as 
\[
\lim_{n\to\infty}[\mu_n(I(n-1))] = \lim_{n\to\infty}\left(\frac{1}{2} - \frac{1}{2^n}\right)
= \frac{1}{2}.
\]
However, no set in  
${\mathcal B}(\overline{\mathbb R})$ is assigned
one half of unit mass by
any of the involved probability measures $\mu, \mu_1, \mu_2, \cdots$.
So $\lim_n[\mu_n(I(n-1))] = 1/2$ is not a probability, and the operation of taking limit is
not closed with respect to probabilities. Nevertheless, $\lim_n[\mu_n(I(n-1))] = 1/2$ is
still the limit of a real sequence. The operation of taking limit remains to be closed with
respect to numbers for any convergent sequence in $[0, 1]$.
$\Box$
}
\end{exam}

\section{Limiting Probabilities on ${\mathcal B}({\mathbb R})$ and
Limits of Real Sequences in $[0, 1]$}
\label{sec-7}
\hskip\parindent
Limiting probabilities are not necessarily limits of probabilities. 
The former must be probabilities. But the latter may or may not be
probabilities. 
A limiting probability is determined by the corresponding limiting 
measure and event. The event may not necessarily be a continuity set of the
limiting measure \cite{Bill2}.
Suppose $\rho_n\Rightarrow\rho$ and $\lim_n\alpha_n = \alpha$
where for all $n$, $\rho_n({\mathbb R}) = 1$ and
$\rho_n(\alpha_n) < 1$ with $\alpha_n\in {\mathcal B}({\mathbb R})$,
$\rho$ is a probability measure on 
${\mathcal B}(\overline{\mathbb R})$ with $\rho(\{-\infty\}) = 0$, and
$\alpha\in {\mathcal B}({\mathbb R})$. If $\lim_n\rho_n(\alpha_n)$
is a probability, $\lim_n[\rho_n(\alpha_n)]$ may not necessarily be a
probability. There is a subtle difference
between $\lim_n\rho_n(\alpha_n)$ and $\lim_n[\rho_n(\alpha_n)]$.
As shown in Example \ref{ex-notight} and Example \ref{ex-limp},
if $\rho({\mathbb R}) = 0$ and if
$\lim_n[\rho_n(\alpha_n)] > 0$ exists, then $\lim_n[\rho_n(\alpha_n)]$ is
not a probability. On the other hand, 
if $\rho({\mathbb R}) = 1$, then $\lim_n[\rho_n(\alpha_n)]$ is
also a probability and coincides with
$\lim_n\rho_n(\alpha_n)$. 
A limiting probability cannot always be expressed as the limit of a real sequence.
Taking limit for a real sequence is not necessarily a closed operation
with respect to probabilities. Nevertheless,
the definition of the limit of a real sequence remains unchanged, 
even if the limit is not a probability. 
We shall prove the above claims in this section.
The case $0< \rho({\mathbb R}) < 1$ may be reduced to
$\rho'({\mathbb R}) = 1$ where $\rho' = \rho/\rho({\mathbb R})$ such that $\rho'_n\Rightarrow\rho'$ and
$\rho'_n = \rho_n/\rho_n({\mathbb R})$ \cite{Loev}. 

To clarify the ambiguity about $\lim_n\rho_n(\alpha_n)$,
we introduce a mapping.
Let $\psi$ be a probability measure 
on ${\mathcal B}(\overline{\mathbb R})$, and
$\eta$ a measurable set in ${\mathcal B}(\overline{\mathbb R})$.
Denote by ${\mathscr P}$ a mapping,
which maps $(\psi, \eta)$ to
$\psi(\eta)$, i.e.,
${\mathscr P}(\psi, \eta)= \psi(\eta)$. 
With ${\mathscr P}$, $\lim_n[\rho_n(\alpha_n)]$ can be expressed by
$\lim_n{\mathscr P}(\rho_n, \alpha_n)$.
If $\lim_n\rho_n(\alpha_n)$ is a probability,
we use $[\lim_n\rho_n(\alpha_n)]_m$ and 
$[\lim_n\rho_n(\alpha_n)]_s$ 
to represent the respective
probability measure and measurable set 
which produce $\lim_n\rho_n(\alpha_n)$.

\begin{rema}
{\em
In general, the argument of $[\cdot]_m$ and
$[\cdot]_s$ corresponds to
a procedure for calculating a probability. The procedure
may simply be a specific probability
assignment such as $\psi(\eta)$, or is given by an expression for computing a
limiting probability, which can eventually be reduced to the
form of $\rho(\alpha)$, where $\rho$ is the measure of the limiting
probability and $\alpha$ the associated $\rho$-measurable set,
determined respectively by weak convergence of the sequence of 
involved probability measures $\rho_n$ and the corresponding measurable 
sets $\alpha_n$ in the context. Procedures for calculating limiting probabilities 
may not necessarily be the same as those for calculating limits of real sequences. 
A procedure for calculating the limit of a real sequence may not necessarily yield 
a limiting probability as shown by Example \ref{ex-limp}. $\Box$
}
\end{rema}

\begin{lemm}
\label{le-1}
Suppose $\rho_n\Rightarrow\rho$ with $\rho(\{-\infty\}) = 0$, $\lim_n\alpha_n = \alpha$, and 
$\rho_n(\alpha_n) < 1$ with $\alpha_n\in {\mathcal B}({\mathbb R})$
and $\rho_n({\mathbb R}) = 1$ for all $n$.

(a) If $\lim_n\rho_n(\alpha_n)$ is a probability, 
\[
\lim_{n\to\infty}\rho_n(\alpha_n) = {\mathscr P}\left(
\rho, \alpha\right).
\]

(b) If both $\lim_n\rho_n(\alpha_n)$ and
$\lim_n{\mathscr P}(\rho_n, \alpha_n)$ are probabilities,
\[
\label{eq-lemb}
\lim_{n\to\infty}\rho_n(\alpha_n) = {\mathscr P}\left(
\rho, \alpha\right)
= \lim_{n\to\infty}{\mathscr P}(\rho_n, \alpha_n).
\]
\end{lemm}

\proof
(a) If $\lim_n\rho_n(\alpha_n)$ is a probability, then
from $\rho_n\Rightarrow\rho$ and $\lim_n\alpha_n=\alpha$
we have
\[
\left[\lim_{n\to\infty}\rho_n(\alpha_n)\right]_m = \rho\;\; 
\mbox{and}\;\;
\left[\lim_{n\to\infty}\rho_n(\alpha_n)\right]_s = \alpha.
\]
Immediately from the definitions of $[\cdot]_m$, $[\cdot]_s$ and ${\mathscr P}$,
\[
\lim_{n\to\infty}\rho_n(\alpha_n) =
\left[\lim_{n\to\infty}\rho_n(\alpha_n)\right]_m
\left(\left[\lim_{n\to\infty}\rho_n(\alpha_n)\right]_s\right)
= \rho(\alpha)={\mathscr P}(\rho, \alpha).
\]

(b) If $\lim_n\rho_n(\alpha_n)$ and
$\lim_n{\mathscr P}(\rho_n, \alpha_n)$ are both probabilities,
we also have $[\lim_n{\mathscr P}(\rho_n, \alpha_n)]_m = \rho$, 
$[\lim_n{\mathscr P}(\rho_n, \alpha_n)]_s = \alpha$, and
\[
\lim_{n\to\infty}{\mathscr P}(\rho_n, \alpha_n) 
= \left[\lim_{n\to\infty}{\mathscr P}(\rho_n, \alpha_n)\right]_m
\left(\left[\lim_{n\to\infty}{\mathscr P}(\rho_n, \alpha_n)\right]_s\right)
= \rho(\alpha)={\mathscr P}(\rho, \alpha)
\]
which together with (a) proves (b).
$\Box$\vspace{10pt}

Thus, if $\lim_n\rho_n(\alpha_n)$ is a probability,
a general expression, or so to speak, a synonym of
$\lim_n\rho_n(\alpha_n)$ is ${\mathscr P}(\rho, \alpha)$.
\begin{rema}
\label{re}
{\em 
By Lemma \ref{le-1}, given $\rho_n\Rightarrow\rho$ and $\lim_n\alpha_n = \alpha$,
if both $\lim_n\rho_n(\alpha_n)$ and $\lim_n[\rho_n(\alpha_n)]$ are probabilities,
then they can be expressed respectively as ${\mathscr P}(\rho, \alpha)$ and
$\lim_n{\mathscr P}(\rho_n, \alpha_n)$, and the order of applying ${\mathscr P}$ and
taking limit is interchangeable.
Labeled by a single index variable $n$, $\lim_n\rho_n(\alpha_n)$ is not a
double limit. A double limit involves a double sequence labeled by two independent
index variables (e.g., $(a_{ij})_{i\geq 1, j\geq 1}$). Issues about
a double limit, such as whether different orders of taking limit will
yield the same result, are irrelevant to $\lim_n\rho_n(\alpha_n)$. $\Box$
}
\end{rema}

The following is a variant of Lemma \ref{le-1}
in some specific settings
with $\alpha_n = I(x_n) = (-\infty, x_n]$, where
$(x_n)_{n\geq 1}$ is a real
sequence, such that 
$\rho_n(I(x_n)) < 1$
for all $n$,
and as $n\to\infty$,
$x_n\uparrow\infty$ or $x_n\to x \in{\mathbb R}$.

\begin{lemm}
\label{le-2}
Assume that $\lim_n\rho_n(I(x_n))$ is a probability,
$\rho_n\Rightarrow\rho$ with $\rho(\{-\infty\}) = 0$, and 
$\rho_n(I(x_n)) < 1$ 
for all $n$.

(a) If $x_n\to x \in{\mathbb R}$ as $n\to\infty$ where $\rho(\{x\}) = 0$, 
then $\lim_n{\mathscr P}(\rho_n, I(x_n))$ is a probability, and
\[
\lim_{n\to\infty}\rho_n(I(x_n)) 
= {\mathscr P}(\rho, I(x))
= \lim_{n\to\infty}{\mathscr P}(\rho_n, I(x_n))
\]
where $I(x) = \lim_nI(x_n) = (-\infty, x]$.

(b) Suppose
$x_n\uparrow\infty$ as $n\to\infty$. Then
\[
\lim_{n\to\infty}\rho_n(I(x_n)) 
= {\mathscr P}(\rho, {\mathbb R}) \leq 1.
\]
If $(\rho_n)_{n\geq 1}$ is tight on
$({\mathbb R}, {\mathcal B}(\mathbb R))$ and if
$\lim_n{\mathscr P}(\rho_n, I(x_n))$ is a probability,
\[
\lim_{n\to\infty}{\mathscr P}(\rho_n, I(x_n)) 
= {\mathscr P}(\rho, {\mathbb R}) = 1.
\]
If $(\rho_n)_{n\geq 1}$ is not tight on
$({\mathbb R}, {\mathcal B}(\mathbb R))$, then
$\lim_n{\mathscr P}(\rho_n, I(x_n))$ 
(if the limit exists) is not a probability, and
\[
{\mathscr P}(\rho, {\mathbb R})
= 1-\rho(\{\infty\})<1.
\]
\end{lemm}

\proof
(a) Since $\rho_n\Rightarrow \rho$ and since $x_n\to x\in{\mathbb R}$,
given $\epsilon > 0$, there exists an $N_{\epsilon}$
such that for all $n \geq N_{\epsilon}$,
\[
|\rho_n(I(x_n)) - \rho(I(x))|\leq 
|\rho_n(I(x_n)) - \rho_n(I(x))| + |\rho_n(I(x)) - \rho(I(x))|
< \epsilon/2 + \epsilon/2 = \epsilon.
\]
So
\[
\lim_{n\to\infty}{\mathscr P}(\rho_n, I(x_n))
= \rho(I(x)) 
= {\mathscr P}(\rho, I(x)).
\]
On the other hand,
\[
\lim_{n\to\infty}\rho_n(I(x_n))
=\left[\lim_{n\to\infty}\rho_n(I(x_n))\right]_m
\left(\left[\lim_{n\to\infty}\rho_n(I(x_n))\right]_s\right)
= \rho(I(x)).
\]
This completes the proof of (a).

(b) Since $\lim_nI(x_n)= \cup_{n\geq 1}I(x_n) = {\mathbb R}$, 
\[
\lim_{n\to\infty}\rho_n(I(x_n)) 
=\left[\lim_{n\to\infty}\rho_n(I(x_n))\right]_m
\left(\left[\lim_{n\to\infty}\rho_n(I(x_n))\right]_s\right)
=\rho({\mathbb R}) = {\mathscr P}(\rho, {\mathbb R}) \leq 1.
\]
If $(\rho_n)_{n\geq 1}$ is tight on
$({\mathbb R}, {\mathcal B}({\mathbb R}))$, 
no mass ``escapes to infinity'' in probability
assignment given by $\rho$, and hence
${\mathscr P}(\rho, {\mathbb R}) =1$.
Moreover,
for each $\epsilon > 0$, there is an interval 
$I_{b(\epsilon)} = (-\infty, b(\epsilon)]$,
such that 
\[
\sup_{n\geq 1}\rho_n\left(I^c_{b(\epsilon)}\right)
= \sup_{n\geq 1}{\mathscr P}\left(\rho_n, I^c_{b(\epsilon)}\right)
< \epsilon
\]
where $I^c_{b(\epsilon)} = {\mathbb R}\setminus I_{b(\epsilon)}$.
For all $n\geq N_{\epsilon}$ where $N_{\epsilon} = \min\{n\in {\mathbb N}: x_n > b(\epsilon)\}$,
$I(x_n) \supset I_{b(\epsilon)}$, 
$I^c(x_n)(= {\mathbb R}\setminus I(x_n))\subset I^c_{b(\epsilon)}$, and hence
\[
{\mathscr P}(\rho_n, I^c(x_n))\leq 
{\mathscr P}\left(\rho_n, I^c_{b(\epsilon)}\right) 
\leq \sup_{n\geq 1}{\mathscr P}\left(\rho_n, I^c_{b(\epsilon)}\right)
< \epsilon.
\] 
In other words, ${\mathscr P}(\rho_n, I^c(x_n))\to 0$ as $n\to\infty$. Consequently,
${\mathscr P}(\rho_n, I(x_n))\to 1$.

Now assume that $(\rho_n)_{n\geq 1}$ is not tight on
$({\mathbb R}, {\mathcal B}({\mathbb R}))$.
In probability assignment given
by $\rho_n$ for each $n$, 
mass must be conserved on the
extended real line. Consequently
\begin{equation}
\label{eq-lem2c}
\rho_n(I(x_n)) 
= \rho_n(\overline{\mathbb R})
- \rho_n(I^c(x_n)) - \rho_n(\{\infty\}).
\end{equation}
Note $\rho_n(I^c(x_n))+\rho_n(\{\infty\})
=\rho_n(\overline{\mathbb R}\setminus I(x_n))$.
In (\ref{eq-lem2c}), $\rho_n(\{\infty\})=0$ and
$\rho_n(\overline{\mathbb R})=1$,
since by assumption, for each $n$, $\rho_n({\mathbb R}) = 1$.
Thus, with ${\mathscr P}$, we can write (\ref{eq-lem2c}) as 
\[
{\mathscr P}(\rho_n, I(x_n)) 
= 1- {\mathscr P}(\rho_n, I^c(x_n)).
\]
Taking limit on both sides of this equation yields
(if the limits exist)
\begin{equation}
\label{eq-lem2clim}
\lim_{n\to\infty}{\mathscr P}(\rho_n, I(x_n))
= 1
- \lim_{n\to\infty}{\mathscr P}(\rho_n, I^c(x_n)).
\end{equation}
If $\lim_n{\mathscr P}(\rho_n, I(x_n))$ is a probability, so
is $\lim_n{\mathscr P}(\rho_n, I^c(x_n))$, such that
\[
\left[\lim_{n\to\infty}{\mathscr P}(\rho_n, I(x_n))\right]_s 
= \lim_{n\to\infty}I(x_n) = {\mathbb R}
\]
\[
\left[\lim_{n\to\infty}{\mathscr P}(\rho_n, I^c(x_n))\right]_s
=\lim_{n\to\infty}I^c(x_n) = \bigcap_{n\geq 1}I^c(x_n) = \emptyset
\]
and (\ref{eq-lem2clim}) implies 
\[
\left[\lim_{n\to\infty}{\mathscr P}(\rho_n, I(x_n))\right]_m
=\left[\lim_{n\to\infty}{\mathscr P}(\rho_n, I^c(x_n))\right]_m
=\rho
\]
with $\rho({\mathbb R}\cup\emptyset)= \rho({\mathbb R})=1$.
However, given that $(\rho_n)_{n\geq 1}$
is not tight on $({\mathbb R}, {\mathcal B}({\mathbb R}))$, we have
$\rho(\{\infty\}) > 0$.
Consequently, $\rho({\mathbb R}) = 1$ implies
\[
\rho({\mathbb R}) + \rho(\{\infty\})
= \rho(\overline{\mathbb R})> 1
\]
and hence violates the 
constraint imposed by the
axioms of probability. As a probability measure on
${\mathcal B}(\overline{\mathbb R})$,
$\rho$ assigns unit mass to $\overline{\mathbb R}$.
So $\lim_n{\mathscr P}(\rho_n, I(x_n))$ is
not a probability. 
Similarly, from
mass conservation on $\overline{\mathbb R}$
in probability assignment given
by $\rho$,
\[
\rho({\mathbb R}) = \rho(\overline{\mathbb R})
- \rho(\{\infty\})
= {\mathscr P}(\rho, {\mathbb R}) = 1 -\rho(\{\infty\})<1
\]
follows immediately. 
$\Box$\vspace{10pt}

\section{The Inconsistency}
\label{sec-8}
\hskip\parindent
Now we demonstrate the inconsistency of measure-theoretic probability.
Let $\overline{\mathbb R}$ be equipped with a topology 
induced by a fixed metric (e.g., see Example \ref{ex-homeo}). 
Let $\lambda, \mu, \lambda_n, \mu_n, n = 1, 2, \cdots$ 
be probability measures on 
$(\overline{{\mathbb R}}, {\mathcal B}(\overline{\mathbb R}))$, 
such that $\lambda, \lambda_1, \lambda_2, \cdots$ are
distributions of i.i.d. random variables $X, X_1, X_2, \cdots$, 
respectively,
as given in Example \ref{ex-tight}, and
$\mu_n$ are distributions of
\begin{equation}
\label{eq-Zn}
Z_n = \max\{i: X_i = Y_n,\; i\in \{1, 2, \cdots, n\}\}
\end{equation}
where
\begin{equation}
\label{eq-Yn}
Y_n = \max\{X_i: i\in\{1, 2, \cdots, n\}\}.
\end{equation}
By (\ref{eq-Zn}), $Z_n, n = 1, 2, \cdots$ form a non-decreasing random sequence $(Z_n)_{n\geq 1}$.
As we can readily see, 
for each $n$, the probability mass function of $Z_n$ is given
by (\ref{eq-1}) in Example \ref{ex-notight}.
As $n\to\infty$, $Z_n$ converges in distribution to an extended random variable
$Z$. Equivalently,
$\mu_n\Rightarrow\mu$, i.e.,
$\mu_n$ converges weakly to $\mu$ induced by $Z$. As a probability measure on
${\mathcal B}(\overline{\mathbb R})$, 
$\mu$ is the Dirac measure concentrated at $\{\infty\}$, i.e.,
$\mu(\{\infty\}) = 1$. 
However, when restricted to ${\mathcal B}({\mathbb R})$, $\mu$ 
is a subprobability measure 
with $\mu({\mathbb R}) = 0$. As shown in Example \ref{ex-tight} and
Example \ref{ex-notight},
$(\lambda_n)_{n\geq 1}$ is tight on 
$({\mathbb R}, {\mathcal B}({\mathbb R}))$ but
$(\mu_n)_{n\geq 1}$ is not (see Sections \ref{sec-4} and \ref{sec-5}).
By (\ref{eq-Zn}) and (\ref{eq-Yn}),
\begin{equation}
\label{eq-Vnn}
Z_n \leq n
\end{equation}
and
\begin{equation}
\label{eq-XVY}
X_{Z_n} = Y_n.
\end{equation}

For the random sequences $(X_n)_{n\geq 1}$ and $(Y_n)_{n\geq 1}$,
a probability measure ${\mathbb P}$ can be defined on the measurable space
$({\mathbb R}^T, {\mathcal B}({\mathbb R}^T))$
where $T = \{1, 2, \cdots\}$.
This probability measure ${\mathbb P}$ can also express probabilities given by the
distribution $\mu_n$ induced by $Z_n$ for each $n$. However, as shown in Section \ref{sec-3},
no probability measure can be defined on $({\mathbb R}^T, {\mathcal B}({\mathbb R}^T))$
to capture the impact of the infinite limit of
$(Z_n)_{n\geq 1}$ on limiting probabilities given by $\mu$. 

In contrast, for the notion of weak convergence, 
probability spaces  
which are the domains
of the involved random variables 
can all be distinct. The domain spaces
are not essential, and actually remain
offstage \cite{Bill2,Loev}.
So we can focus on the probability measures
and their weak limits, as long as the range space and the topology on it are the same for
all the random variables involved 
\cite{Bill2}.
The weak limit of $\lambda_n$ is trivially $\lambda$, i.e., 
$\lambda_n\Rightarrow\lambda$
with $\lambda(\{0\}) = 1/2$.
However, the following theorem is also provable.

\begin{theo}
\label{th}
For $\lambda_n$ given in Example \ref{ex-tight}, $\lambda_n\Rightarrow\lambda$ with 
$\lambda(\{0\}) = \mu({\mathbb R}) = 0$. 
\end{theo}

\proof
By (\ref{eq-Vnn}), 
\begin{equation}
\label{eq-2p}
\lambda_n(\{0\}) = {\mathbb P}(X_n = 0) =  
{\mathbb P}(X_n = 0, Z_n < n) +
{\mathbb P}(X_n = 0, Z_n = n).
\end{equation}
Consider the limiting probabilities on both
sides of (\ref{eq-2p}).
We first show $\lim_n{\mathbb P}(X_n =0, Z_n = n)=0$.
\begin{eqnarray}
\label{eq-rec}
&&\lim_{n\rightarrow\infty}
{\mathbb P}(X_n =0, Z_n = n) \nonumber\\
&& =
\lim_{n\rightarrow\infty}
{\mathbb P}(X_{Z_n} =0, Z_n = n)\nonumber \\
&& =
\lim_{n\rightarrow\infty}
{\mathbb P}(Y_n =0, Z_n = n)
\;(\mbox{by (\ref{eq-XVY}))} \nonumber\\
&& = 0. 
\end{eqnarray}
The last equality follows from 
${\mathbb P}(\lim_{n\to\infty}Y_n = 1)= 1$,
\[
\{Y_n =0, Z_n = n\}
\subset \{Y_n =0\},\;\;
{\mathbb P}(Y_n =0, Z_n = n)
\leq {\mathbb P}(Y_n = 0)
\]
and 
\[
\lim_{n\to\infty}{\mathbb P}(Y_n= 0)
= {\mathbb P}\left (\bigcap_{n\geq 1}\{Y_n=0\}\right )
(\mbox{since $\{Y_n =0\}\supset\{Y_{n+1}=0\}$
as implied by (\ref{eq-Yn})}).
\]
Similarly, 
from (\ref{eq-2p}), (\ref{eq-rec}), and 
\[
\{X_n = 0, Z_n < n\}
= \{Z_n < n\}
\]
we have 
\[
{\mathbb P}\{X_n =0, Z_n < n\} = {\mathbb P}\{Z_n < n\}
=\mu_n(I(n-1))
\]
and 
\[
\lambda(\{0\}) = \lim_{n\to\infty}\lambda_n(\{0\})
= \lim_{n\to\infty}{\mathbb P}(X_n=0, Z_n < n)
= \lim_{n\to\infty}\mu_n(I(n-1)). 
\]
Since $\lambda(\{0\})$ is a probability, so is $\lim_n\mu_n(I(n-1))$, and
\[
\lambda(\{0\}) 
=\lim_{n\to\infty}\mu_n(I(n-1)) 
= {\mathscr P}(\mu, {\mathbb R}) = \mu({\mathbb R}) = 0
\]
follows from
Lemma \ref{le-2}(b) with $x_n = n-1$, 
$\rho_n = \mu_n$ and $\rho(\{\infty\})= \mu(\{\infty\})=1$.
$\Box$\vspace{10pt}
\begin{rema}
\label{re-2}
{\em
By definition, ``1/2'' in 
the probability mass function ``$\lambda(\{0\})= 1/2$'' (see Example \ref{ex-tight})
represents a (numerical) probability. In contrast, 
as explained in Example \ref{ex-notight} and Example \ref{ex-limp},
``1/2'' in ``$\lim_n[\mu(I(n-1))]= \lim_n{\mathscr P}(\mu_n, I(n-1)) = 1/2$'' is not
a probability. As we have discussed in Sections \ref{sec-6} and \ref{sec-7},
without representing the value of a quantity in the physical world, a real number
is only a point in a set formulated by a system of axioms. It makes no sense at all
to treat a probability, which is a quantified attribute of an event in the
physical world, as a point
of a set, which is merely a mathematical 
(or metaphysical)
entity still considered questionable 
by some mathematicians even nowadays \cite{Brid,Rose}. 
Since  $\lambda(\{0\}) = \lim_n\mu(I(n-1))$ as shown in the proof of
Theorem \ref{th}, and since $\lambda(\{0\})$ is a probability, 
$\lim_n\mu(I(n-1))$ must be a probability, and this probability is
${\mathscr P}(\mu, {\mathbb R}) = \mu(\mathbb R) = 0$ by Lemma \ref{le-2}. Since ``1/2'' in
``$\lim_n[\mu(I(n-1))]= \lim_n{\mathscr P}(\mu_n, I(n-1)) = 1/2$''
is a real number, i.e., a point on the abstract real line ${\mathbb R}$ but
not a probability, 
it is meaningless to equate $\lambda(\{0\})$ with $\lim_n[\mu(I(n-1))] = \lim_n{\mathscr P}(\mu_n, I(n-1))$
(see also Example \ref{ex-unitpro} and Example \ref{ex-info}). 
$\Box$
}
\end{rema}

\section{Concluding Remarks}
\label{sec-9}
\hskip\parindent
In general, the presupposition 
for questions of identity (expressed by ``$=$'') in a theory is that the two sides
inquired about both belong to the same category \cite{Bena2}.
An expression ``$x = y$'' where
$x$ and $y$ are in the same category but numerically different is erroneous 
only when the theory is proved to be contradiction-free as requested by Hilbert.
Unfortunately, according to G{\" o}del's theorems, 
measure-theoretic probability is one of the theories for which such a proof
does not exist. Nevertheless, if one side of ``$=$'' is a probability, then
the other side must also be a probability. This requirement is necessary
to prevent confusion of probabilities and arbitrarily calculated numbers
which are not eligible to serve as probabilities. 
Whether a number is eligible to serve as a probability is 
determined by whether the number is the result of 
applying a probability measure to a measurable set in the context.
In particular, when we calculate
limiting probabilities, the two sides of ``$=$'' must both be limiting probabilities.
Since procedures for calculating limits of real sequences 
may not necessarily yield limiting probabilities, the limit of a real sequence which
is not a probability is irrelevant to whatever implied by calculations of
limiting probabilities.

Consider 
$\lim_n\mu_n(I(n-1))$ (see the proof of Theorem \ref{th}). 
The probability measures $\mu_n$ and the associated measurable sets 
$I(n-1)$ approach their respective limits simultaneously. 
This implies two different orders of operations for calculating
$\lim_n\mu_n(I(n-1))$:
\begin{itemize}
\item[(a)]
taking the limits of the probability measures $\mu_n$
and the measurable sets $I(n-1)$ 
simultaneously first, then applying the limiting measure $\mu$ to 
the limiting measurable set ${\mathbb R}$, i.e.,
\[
\lim_{n\to\infty}\mu_n(I(n-1)) = {\mathscr P}(\mu, {\mathbb R}) = \mu({\mathbb R})
\]
\item[(b)] 
applying the measures $\mu_n$ to the corresponding measurable sets 
$I(n-1)$ to obtain a sequence $(\mu_n(I(n-1)))_{n\geq 1}$ of 
numbers first, then taking the limit of the real sequence, i.e.,
\[
\lim_{n\to\infty}\mu_n(I(n-1))=\lim_{n\to\infty}{\mathscr P}(\mu_n, I(n-1)).
\]
\end{itemize}
Order (a) produces a probability, but order 
(b) yields a number which is not a probability (see Example \ref{ex-limp}). So the meaning of 
$\lim_n\mu_n(I(n-1))$ is ambiguous. 

Since $\lim_n\mu_n(I(n-1))$ can only be written as a ``linear'' string of symbols,
such a ``linear'' way of writing might result in an impression as if 
(b) were the only option for calculating $\lim_n\mu_n(I(n-1))$. 
However, when standing in the relation ``$=$'' to
a probability (i.e., $\lambda(\{0\})$), $\lim_n\mu_n(I(n-1))$ has to be calculated according to
(a) rather than (b) to yield a probability as well, for
the two sides of ``$=$'' must both be
probabilities, even if they are numerically different. 
As we have
demonstrated, ``two numerically different
probabilities equal each other'' is a provable statement, which shows that
measure-theoretic probability is inconsistent. 

The revealed contradiction of measure-theoretic probability does not exist in the physical world.
It appears only in the theory. Practical applications of
measure-theoretic probability will not be affected 
by the inconsistency as long as ideal events in the theory
are not mistaken for real events in the physical world. On the other hand,
the inconsistency of measure-theoretic probability cannot be explained away, and must be resolved.
Constructive mathematics (e.g., see \cite{Brid,Rose}) can avoid such inconsistency. There is no contradiction 
reported in constructive mathematics.

\end{document}